\documentclass{amsart}
\usepackage{amscd}
\usepackage{amsmath}
\usepackage{amsxtra}
\usepackage{amsfonts}
\usepackage{amssymb}

\newtheorem{theorem}{Theorem}[section]
\newtheorem{corollary}[theorem]{Corollary}
\newtheorem{lemma}[theorem]{Lemma}
\newtheorem{proposition}[theorem]{Proposition}
\theoremstyle{definition}
\newtheorem{definition}[theorem]{Definition}
\newtheorem{remark}[theorem]{Remark}

\newtheorem{example}[theorem]{Example}
\theoremstyle{remark}

\renewcommand{\theclaim}{\textup{\theclaim}}

\newtheorem*{acknowledgements}{Acknowledgements}

\numberwithin{equation}{section}

\def\openone

{\mathchoice

{\hbox{\upshape \small1\kern-3.3pt\normalsize1}}

{\hbox{\upshape \small1\kern-3.3pt\normalsize1}}

{\hbox{\upshape \tiny1\kern-2.3pt\SMALL1}}

{\hbox{\upshape \Tiny1\kern-2pt\tiny1}}}

\makeatletter

\newbox\ipbox

\newcommand{\ip}[2]{\left\langle #1\,|\,#2\right\rangle}

\newcommand{\diracb}[1]{\left\langle #1\mathrel{\mathchoice

{\setbox\ipbox=\hbox{$\displaystyle \left\langle\mathstrut
#1\right.$}

\vrule height\ht\ipbox width0.25pt depth\dp\ipbox}

{\setbox\ipbox=\hbox{$\textstyle \left\langle\mathstrut
#1\right.$}

\vrule height\ht\ipbox width0.25pt depth\dp\ipbox}

{\setbox\ipbox=\hbox{$\scriptstyle \left\langle\mathstrut
#1\right.$}

\vrule height\ht\ipbox width0.25pt depth\dp\ipbox}

{\setbox\ipbox=\hbox{$\scriptscriptstyle \left\langle\mathstrut
#1\right.$}

\vrule height\ht\ipbox width0.25pt depth\dp\ipbox}

}\right. }

\newcommand{\dirack}[1]{\left. \mathrel{\mathchoice

{\setbox\ipbox=\hbox{$\displaystyle \left.\mathstrut
#1\right\rangle$}

\vrule height\ht\ipbox width0.25pt depth\dp\ipbox}

{\setbox\ipbox=\hbox{$\textstyle \left.\mathstrut
#1\right\rangle$}

\vrule height\ht\ipbox width0.25pt depth\dp\ipbox}

{\setbox\ipbox=\hbox{$\scriptstyle \left.\mathstrut
#1\right\rangle$}

\vrule height\ht\ipbox width0.25pt depth\dp\ipbox}

{\setbox\ipbox=\hbox{$\scriptscriptstyle \left.\mathstrut
#1\right\rangle$}

\vrule height\ht\ipbox width0.25pt depth\dp\ipbox}

} #1\right\rangle}

\newcommand{\cj}[1]{\overline{#1}}

\newcommand{\bz}{\mathbb{Z}}
\newcommand{\br}{\mathbb{R}}
\newcommand{\bc}{\mathbb{C}}
\newcommand{\bt}{\mathbb{T}}
\newcommand{\bn}{\mathbb{N}}

\newcommand{\cc}{\mathfrak{c}}
\usepackage{graphicx}

\def\blfootnote{\xdef\@thefnmark{}\@footnotetext}

\begin{document}
\title[Hilbert spaces, similarity, dynamics]{Hilbert spaces built on a similarity and on dynamical renormalization}
\author{Dorin Ervin Dutkay}
\address[Dorin Ervin Dutkay]{Department of Mathematics\\
Rutgers, The State University of New Jersey\\
Hill Center-Busch Campus\\
110 Frelinghuysen Road\\
Piscataway, NJ 08854}
\email{ddutkay@math.rutgers.edu}\ 

\author{Palle E.T. Jorgensen}
\address[Palle E.T. Jorgensen]{Department of Mathematics\\
The University of Iowa\\
14 MacLean Hall\\
Iowa City, IA 52242}
\email{jorgen@math.uiowa.edu}

\thanks{Work supported in part by the U.S. National Science
Foundation under grants DMS-0457491 and DMS-0457581}
\begin{abstract}
 We develop a Hilbert-space framework for a number of general multi-scale problems from dynamics. The aim is to identify
 a spectral theory for a class of systems based on iterations of a non-invertible endomorphism.
 We are motivated by the more familiar approach to wavelet theory which starts with the two-to-one endomorphism
 $r\colon z \mapsto z^2$ in the one-torus $\bt$, a wavelet filter, and an associated transfer operator. This leads to a
 scaling function and a corresponding closed subspace $V_0$ in the Hilbert space $L^2(\br)$.  Using the dyadic
 scaling on the line $\br$, one has a nested family of closed subspaces $V_n$, $n \in \bz$, with trivial intersection,
 and with dense union in $L^2(\br)$. More generally, we achieve the same outcome, but in different Hilbert spaces,
 for a class of non-linear problems. In fact, we see that the geometry of scales of subspaces in Hilbert space is
 ubiquitous in the analysis of multiscale problems, e.g., martingales, complex iteration dynamical systems,
 graph-iterated function systems of affine type, and subshifts in symbolic dynamics. We develop a general
 framework for these examples which starts with a fixed endomorphism $r$ (i.e., generalizing  $r(z) = z^2$ )
 in a compact metric space $X$. It is assumed that $r \colon X\rightarrow X$ is onto, and finite-to-one.
\end{abstract}
\maketitle \tableofcontents
\section{Introduction}
We study a class of endomorphisms $r \colon X\rightarrow X$, where $X$ is
a metric space. The endomorphism is assumed onto, and finite-to-one.
We build a spectral theory on a Hilbert space associated naturally
with $(X, r)$. Our focus is on the case when $X$ is assumed to carry
a certain strongly invariant measure $\rho$, see (\ref{eqstinv}).
\par
Continuing our earlier work \cite{DuJo05} we consider
basis constructions in a general context of dynamical systems; the case
of endomorphisms, i.e., non-reversible dynamics. Our framework will
include wavelet bases, as well as algorithmic basis constructions in
Hilbert spaces built on fractals or on Julia sets of rational functions
in one complex variable. In fact, these examples motivated our results.

          First recall that in the real variable case of standard
wavelets (in one or several variables, i.e., the $d$-dimensional Lebesgue
measure), there is a separate generalizations of standard dyadic
wavelets, again based on translation and scaling: See for example
\cite{BJMP05} for such an approach to the construction of generalized wavelet
bases in the Hilbert space $L^2(\br^d)$, i.e., of orthogonal bases in $L^2(\br^d)$, or just frame wavelet bases, but still in $L^2(\br^d)$.

It is the purpose of this paper to develop a geometric context of this
viewpoint which applies to any kind of dynamics which is based on an
iterated scale of selfsimilarity. Hence our paper will offer a Hilbert-%
space framework which goes beyond the setting of scale similarity, and
our results will offer a new viewpoint even in the case of the more
familiar selfsimilarity which is based on a cascade of affine scales.

         The best know instance of this is $d = 1$, and dyadic wavelets
\cite{Dau92}. In that case, the two operations on the real line $\br$ are
translation by the group $\bz$ of the integers, and scaling by powers of 2,
i.e., $x\mapsto 2^jx$, as $j$ runs over $\bz$.  This is the approach to wavelet theory
which is based on multiresolutions and filters from signal processing. In
higher dimensions $d$, the scaling is by a fixed matrix, and the
translations by the rank-$d$ lattice $\bz^d$. Again we will need scaling by all
integral powers. We view points $x$ in $\br^d$ as column vectors, and we then
consider the group of scaling transformations, $x\mapsto A^j x$  as $j$ ranges
over $\bz$.

 Suitable spectral conditions will be imposed on $A$. In particular
we note that if $A$ is integral, i.e., the entries in $A$ are in $\bz$, then $x\mapsto A x$ passes to the quotient $\br^d/\bz^d$. Since $\br^d/\bz^d$ is a copy of the
compact $d$-torus $\bt^d$ via a familiar identification, we see that $A$ induces
an endomorphism $r_A$ in $\bt^d$. If further $A$ is invertible, then $r_A$ is
finite-to-one, and maps $\bt^d$ onto itself. In fact, for every $x$ in $\bt^d$, the
inverse image $r_A^{-1}(x)$ has cardinality $= | \det A | =: N$.

\begin{equation}\label{eq1_1}
\frac{1}{\sqrt{\det A}}\varphi(A^{-1}x)=\sum_{k\in\bz^d}a_k\varphi(x-k),\quad(x\in\br^d).
\end{equation}
 
So our starting point is a given finite-to-one endomorphism $r\colon X \rightarrow X$ in a
compact space $X$. Our aim is three-fold: (1) to build an associated
Hilbert space which admits wavelet decompositions; (2) to show that the
corresponding computations can be done with a geometric algorithm; and
finally (3) we offer concrete examples from dynamics where our approach
leads to new insight. So in addition to the endomorphism $(X, r)$, our
initial setup will include a scalar function $m_0$; an analogue of the
function from wavelet theory which determines low-pass filters.

     Details: Set $W(x):= |m_0 |^2 / \# r^{-1}(x) $. We say that $m_0$
satisfies a {\it low-pass} condition if $W(0) = 1$.
(In the special case of (\ref{eq1_1}) above, the relationship between the
function $m_0$ and the coefficients $\{a_k\}$ is that the $a_k$ numbers will be
the $d$-Fourier coefficients of $m_0$ when $m_0$ is viewed as a function on the
compact quotient $X = \br^d/\bz^d$. This explains the summation over $\bz^d$ in
(\ref{eq1_1}).)

Suppose for some $p$, and $x \in X$, that $r^p(x) = x$. Then we say that the
finite set of points $C = \{ x, r(x),\dots, r^{p-1}(x) \}$ is a {\it cycle}. A
cycle $C$ is called a {\it $W$-cycle} if $W(y) = 1$ for all $y\in C$.

 We will extend to the context of endomorphisms the following general
principle from wavelets in the Hilbert space $L^2(\br^d)$: A generalized
wavelet basis (also called a Parseval-frame, see e.g., \cite{BJMP05}) will have
the stronger orthonormal basis (ONB)property when the only one $W$-cycle is
$C = \{0\}$. On the other hand, the presence of non-trivial $W$-cycles is
consistent with wavelet systems that form frame-bases. The reader is
referred to \cite{BDP05} for details regarding these more general
wavelet bases. It was proved in \cite{BDP05} that the
presence of $W$-cycles is consistent with a class of certain
{\it super-wavelets}. This wavelet basis involves an additional cyclic
structure which we will develop in the paper.

\par
This setup arose earlier for the familiar linear multiresolution
analysis (MRA) approach to wavelets: Recall \cite{Dau92} that dyadic
wavelets represent a special basis for the Hilbert space $L^2(\br)$,
but they are generated by a subspace $V_0$ in $L^2(\br)$ which is
the closed linear span of a single function $\varphi$ and its
translates by the integers $\bz$. The function $\varphi$ satisfies a
certain scaling identity
\begin{equation}\label{eqscalr}
\frac{1}{\sqrt{2}}\varphi(x/2)=\sum_{k\in\bz}a_k\varphi(x-k),\quad(x\in\br).
\end{equation}
which implies that the scaling operator $U f(x) := 1/\sqrt{2}
f(x/2)$ maps $V_0$ into itself. A solution $\varphi$ is called a
scaling function. Using a terminology from optics, we say that
functions on $\br$ represent signals or images, and that the
subspace $V_0$ initializes a fixed resolution.
\par
\begin{figure}[h]
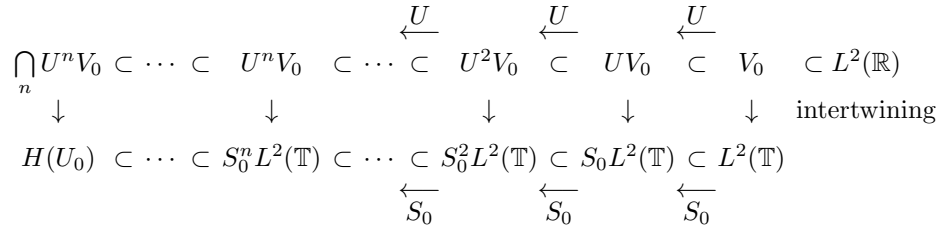

$\setlength{\arraycolsep}{2pt}
\begin{array}{cccccccccccccl}
&&&&&&& \makebox[0pt]{\hss$\overset{\textstyle
U}{\longleftarrow}$\hss} && \makebox[0pt]{\hss$\overset{\textstyle
U}{\longleftarrow}$\hss}
&& \makebox[0pt]{\hss$\overset{\textstyle U}{\longleftarrow}$\hss} && \\
\smash{\bigcap\limits_n U^n V_0} & \subset & \cdots & \subset & U^n V_0 & \subset & \cdots & \subset & U^2 V_0 & \subset & U V_0 & \subset & V_0 & {} \subset L^2(\mathbb{R}) \\[2\jot]
 \downarrow &&&& \downarrow &&&& \downarrow && \downarrow && \downarrow & \text{intertwining} \\[2\jot]
H(U_0) & \subset & \cdots & \subset & S_0^n L^2(\mathbb{T}) & \subset & \cdots & \subset & S_0^2 L^2(\mathbb{T}) & \subset & S_0 L^2(\mathbb{T}) & \subset & L^2(\mathbb{T}) & \\
&&&&&&& \makebox[0pt]{\hss$\underset{\textstyle
S_0}{\longleftarrow}$\hss} &&
\makebox[0pt]{\hss$\underset{\textstyle S_0}{\longleftarrow}$\hss}
&& \makebox[0pt]{\hss$\underset{\textstyle
S_0}{\longleftarrow}$\hss} &&
\end{array}
$ \caption{Multiresolutions} \label{FigMultiresolutions}
\end{figure}

A special case: $X = \bt = \{z \in \bc\, |\,  |z| = 1\} = \br/\bz$,
$r(z) = z^2$, and $m_0$ is the function on $\bt$  with Fourier
coefficients equal to the masking coefficients $a_k$ from
(\ref{eqscalr}), i.e., $m_0(z)=\sum_{k\in\bz}a_kz^k$. The function
$m_0$ is called a filter function because of an analogy to a setting
in signal processing. One of the axioms for $m_0$ (the {\it
quadrature-mirror-filter} axiom) from wavelet theory amounts to the
fact that the associated linear operator, $S_0 h (z) : = m_0(z)
h(z^2)$  is isometric in $L^2(\bt, \mbox{ Haar measure})$; see
Figure \ref{FigMultiresolutions}.

\par In this paper, we state a version of the scaling identity
(\ref{eqscalr}), for the case of an endomorphism $r\colon X\rightarrow X$,
and we show that it admits a solution in certain Hilbert spaces
built on $(X,r)$. It turns out that the variant of (\ref{eqscalr})
which arises by the Fourier transform, i.e.,
\begin{equation}\label{eqscalir}\sqrt{2}\hat\varphi(t)=m_0(e^{it/2})\hat\varphi(t/2),\quad(t\in\br),
\end{equation}
is more suggestive of the generalization we have in mind; see
Theorem \ref{propphic} for details.

\par
While the standard MRA approach to wavelets (see \cite{HeWe96})
restricts the functions $m_0$ in (\ref{eqscalir}) by assuming that
$m_0$ is in some regularity class, e.g., is Lipschitz, we shall not
do this here. Moreover, there is a rich class of wavelet systems
where $m_0$ is typically only known {\it a priori} to be $L^\infty$.
This is the case, for example, for the frequency localized wavelets
studied in \cite{BJMP05} and \cite{LPT01} (In this last case, $m_0$
is in fact matrix-valued.)
\par The scaling identity (\ref{eqscalr}) implies that there is a
natural intertwining of the isometry $S_0$ on $L^2(\bt)$ with the
restriction of $U$ to the subspace $V_0$ in $L^2(\br)$. A second
axiom for $m_0$ from wavelet theory (called {\it low-pass}) implies
that $S_0$ is a {\it pure} shift isometry, i.e., that the
intersection $S_0^n(L^2(\bt))$, for $n$ in $\bn$ is $\{0\}$, see
Figure \ref{FigMultiresolutions}. Because of the intertwining
relation, this fact guarantees that the standard functions that make
up a wavelet basis really do form a basis for the whole Hilbert
space $L^2(\br)$. See Figure \ref{FigMultiresolutions}. And the
purity of $S_0$ is also what yields a certain martingale system,
i.e., a nested family of spaces, or of sigma-algebras.
\par
It is the purpose of this paper to generalize this setting to that
of endomorphisms, and to realize a natural scaling function, as a
generating vector in a Hilbert space which corresponds to $L^2(\br)$
for the special case of wavelets. For this purpose we introduce a
solenoid $X_\infty$ built on $X$, and a family of repelling cycles
for the system $(X, r, m_0)$. Our Hilbert space is built as an
$L^2$-space on certain infinite paths starting at $X$. In Theorem
\ref{propphic}, we solve the corresponding scaling identity, and
write the scaling function as an infinite product. As one should
expect by analogy to wavelets, a central theme in our present
analysis is a characterization of those filter functions $m_0$ on
$X$ for which the scaling identity has non-trivial solutions in a
Hilbert space of functions of $X_\infty$.
\par
 A concrete example of this wavelet technique used on a particular graph dynamical system (The Golden mean shift) is presented in Proposition \ref{pr2_18} below. Our aim is to present this as a systematic tool for dynamics outside the traditional context of wavelets in $L^2(\br)$.
 \par In recent
papers \cite{DuJo03, DuJo05}, the co-authors have adapted this MRA
technique to a related but different problem, the  problem of
creating a spectral theory for a class of non-linear iterated
function systems (IFS), but in those cases, there is not a direct
analogue to the scaling identity. Our construction here parallels
the one we outlined briefly for the standard dyadic MRA wavelet
constructions \cite{HeWe96}. (We have sketched the standard wavelet
construction only in the dyadic case, and only in one dimension,
i.e., for $\br$, but it is known that this construction carries over
{\it mutatis mutandis} to $\br^d$ with $d > 1$, and when $x\mapsto
2x$, is replaced with matrix scaling $x\mapsto A x$ in $\br^d$ where
$A$ is a $d$ by $d$ matrix over $\bz$ with eigenvalues $\lambda$
such that $| \lambda |> 1$. Moreover our results apply to the kind
of multiwavelets studied recently in \cite{BJMP05}.) \par Our
present paper is not about $\br^d$-wavelets but instead about a
class of non-linear dynamics $r\colon X\rightarrow X$.
 Specifically, now we start with $r\colon X\rightarrow X$, and the function $m_0$ is defined on $X$. We will also call $m_0$
 a {\it filter function} because of known analogy to subband filtering in signal processing. When $m_0$ is given then
 $S_0$ given by  $S_0h(x):= m_0(x) h(r(x))$ is isometric in $L^2(X)$, subject to a technical condition on $m_0$. So by Wold's theorem \cite{NaFo70}, it is then the orthogonal sum of a
 shift operator $S$ and a unitary operator $U_0$; i.e., the Hilbert space $L^2(X)$ on which $S_0$ acts is the direct sum of two
 Hilbert spaces $H(S)$ and $H(U_0)$ such that (i) each space invariant for $S_0$, (ii) the restriction to $H(S)$ is a
 shift $S$, and the restriction to $H(U_0)$ a unitary operator. We say that $S_0$ is {\it pure} if $H(U_0) =
 \{0\}$. (See Figure \ref{FigMultiresolutions}.) This will be equivalent to the fact that the intersection
 of the multiresolution subspaces is trivial.

\par
This means that $S_0$ is itself a shift operator on $L^2(X)$. Our
theorem \ref{thwold} gives a simple condition for the isometry $S_0$
to be pure.
\par
The first step in our construction is an extension from the initial
endomorphism $r\colon X \rightarrow X$, to a new invertible system $\hat
r \colon X_\infty \rightarrow X_\infty$, i.e., with $\hat r$ invertible
on $X_\infty$. When $r$ is assumed finite-to-one, this can be done
such that there is a quotient  $X_\infty/X$ which becomes a Cantor
space. The extension space $X_\infty$ is called a solenoid. For the
case when $(X, r)$ is a one-sided subshift \cite{MaUr03}, we work
out (in Section \ref{subsscal}) an explicit model for this solenoid.
\par
In fact the notion of a solenoid (for the study of dynamics of an
endomorphism and extension to an automorphism) was used already in a
pioneering paper by Lawton \cite{Law73} in 1973. Lawton considered
groups with expansive automorphisms; see also \cite{LaPa67}.
Motivated by applications, we note that our present analysis is not
restricted to groups.
\par
The use of solenoids in the study of particular systems with scale
similarity was initiated in the paper \cite{BrJo91}, and was
continued in \cite{Bre96}. The context of \cite{BrJo91} is a class
of algebraic irrational numbers and an associated $C^*$-algebraic
crossed product. In a general context of non-linear dynamics, this
work was continued in \cite{DuJo03, DuJo05}.

\section{Covariant representations}\label{scova}
\par
Let $X$ be a compact metric space with a non-invertible endomorphism
$r\colon X\rightarrow X$ such that $r$ is measurable, onto and finite to
one, i.e., $0<\#r^{-1}(x)<\infty$ for all $x\in X$.
\par
We have shown in \cite{DuJo04} and \cite{DuJo05} that, for certain
filter functions $m_0$ on $X$, one can construct multiresolutions
and scaling functions in Hilbert spaces of functions on $X_\infty$
(see (\ref{eqxinfty})).\par In \cite{DuJo05} we proved that to get
useful multiresolutions, the function $m_0$ must have certain
extreme cycles (see Definition \ref{def2_10}). In this case the
measure on $X_\infty$ is actually supported on a smaller set
$\mathbf{N}_C$ (see (\ref{eqnc}) below).
\subsection{The ground space}\label{subsgrou}
An (infinite) {\it path} starting at $x$ is a sequence
$(z_1,z_2,\dots)$ of points in $X$ such that $r(z_1)=x$,
$r(z_{n+1})=z_n$ for $n\geq1$. We denote by $\Omega_x$ the set of
paths starting at $x$. We denote by $X_\infty$ the set of all paths,
\begin{equation}\label{eqxinfty}
X_\infty=\cup_{x\in X}\Omega_x.\end{equation}
\par
Note that a path $(z_1,z_2,\dots)$ in $\Omega_x$ can be identified
with the doubly infinite sequence $(z_n)_{n\in\bz}$, where $z_0:=x$
and $z_{-n}=r^n(x)$ for $n\geq0$.
\par
$X_\infty\subset X^{\bz}$ inherits the usual Tychonoff topology from
$X^{\bz}$.
\par
The maps $\theta_n\colon X_\infty\rightarrow X$ are defined for all
$n\in\bz$, by
$$\theta_n((z_k)_{k\in\bz})=z_n.$$
\par
The endomorphism $r$ can be extended to the automorphism $\hat r$
defined on $X_\infty$ by
$$\hat r(z_n)_{n\in\bz}=(z_{n-1})_{n\in\bz}.$$
\par
These maps satisfy the following relations:
$$\theta_{n}\circ\hat r=\theta_{n-1},\quad \theta_0\circ\hat r=r\circ\theta_0.$$
\par
For a function $g$ on $X$ we define
\begin{equation}\label{eqglan}
g^{(n)}(x):=g(x)g(r(x))\cdots g(r^{n-1}(x)),\quad(n\geq1).
\end{equation}
\par
For a function $\xi$ on $X_\infty$, we define $\xi^{(0)}=1$,
$$\xi^{(n)}:=\xi\xi\circ\hat r\cdots\xi\circ\hat r^{n-1},$$
and if $\xi$ is not vanishing on $X_\infty$, then
$$\xi^{(-n)}=\frac{1}{\xi\circ\hat r^{-1}\xi\circ\hat
r^{-2}\cdots\xi\circ\hat r^{-n}},\quad(n\geq1).$$
\par
We can identify functions $g$ on $X$ with functions on $X_\infty$ by
$g\leftrightarrow g\circ\theta_0$. (Note that the two definitions
for $g^{(n)}$ will coincide.)

Consider $r\colon X\rightarrow X$ and suppose $\rho$ is a {\it strongly
invariant} probability measure on $X$, i.e.,
\begin{equation}\label{eqstinv}
\int_Xf(x)\,d\rho(x)=\int_X\frac{1}{\#r^{-1}(x)}\sum_{y\in
r^{-1}(x)}f(y)\,d\rho(x),\quad(f\in L^\infty(\rho)).\end{equation}
\par
Let $C=\{x_0,x_1,\dots,x_{p-1}\}\subset X$ be a {\it cycle} of length
$p$, i.e., the points $x_i$ are distinct and $r(x_{i+1})=x_i$,
$r(x_0)=x_{p-1}$.
\par
We define the set
\begin{equation}\label{eqnc}
\mathbf{N}_C(x):=\{\omega=(z_1,z_2,\dots)\in\Omega_x\,|\, \lim_{n\rightarrow\infty}z_{pn}\in C\}.
\end{equation}
For each $x\in X$ and $\omega=(z_1,z_2,\dots)\in\mathbf{N}_C(x)$, define $i(\omega)\in\{0,\dots,p-1\}$ by $i(\omega):=i$ if $\lim_{k\rightarrow\infty}z_{kp}=x_i$.
\par
An inspection reveals that each $\mathbf{N}_C(x)$ is countable.
\par
Define the measure $\lambda_C$ on $X_\infty$ by
\begin{equation}\label{eqlambdac}
\int_{X_\infty}f\,d\lambda_C=\int_X\sum_{\omega\in\mathbf{N}_C(x)}f(\omega)\,d\rho(x).
\end{equation}

To simplify the notation we write $\cc(x)=\#r^{-1}(r(x))$.
\begin{proposition}\label{proplambda}
For all $\xi\in L^1(X_\infty,\lambda_C)$ and $n\in\bz$, we have
$$\int_{X_\infty}\cc^{(n)}\xi\circ\hat r^n\,d\lambda_C=\int_{X_\infty}\xi\,d\lambda_C.$$
\end{proposition}

\begin{proof}
It is enough to prove this for $n=1$, the rest follows by induction.
$$\int_{X_\infty}\cc\xi\circ\hat r\,d\lambda_C=\int_X\#r^{-1}(r(x))\sum_{\omega\in\mathbf{N}_C(x)}\xi(\hat r(x,\omega))\,d\rho(x)=$$
$$\int_X\frac{1}{\#r^{-1}(x)}\sum_{y\in r^{-1}(x)}\#r^{-1}(r(y))\sum_{\omega\in\mathbf{N}_C(y)}\xi(\hat r(y,\omega))\,d\rho(x)=$$
$$\int_X\sum_{\omega\in\mathbf{N}_C(x)}\xi(x,\omega)\,d\rho(x)=\int_{X_\infty}\xi\,d\lambda_C.$$
\end{proof}

\subsection{The operators}\label{subsoper}
\par
In this subsection we show that when $(X,r)$ is given as above,
there is a natural covariant representation $(U,\pi)$ acting on the
Hilbert space $L^2(X_\infty,\lambda_C)$, i.e., with $r$ inducing a
unitary operator $U$, and $\pi$ a representation of $L^\infty(X)$ by
multiplication operators, such that the relation (\ref{eqcov}) is
satisfied on $L^2(X_\infty,\lambda_C)$. (The operators on
$L^2(X_\infty,\lambda_C)$ are equipped with the strong operator
topology (SOT).)
\par
Let $\alpha_0,\dots,\alpha_{p-1}$ be a set of complex numbers of
absolute value $1$.
\par
Let $U$ be the operator on $L^2(X_\infty,\lambda_C)$ defined by
\begin{equation}\label{equc}
U\xi(x,\omega)=\alpha_{i(\omega)}\sqrt{\#r^{-1}(r(x))}\xi\circ\hat
r(x,\omega),\quad(\xi\in L^2(X_\infty,\lambda_C),x\in
X,\omega\in\Omega_x).
\end{equation}
For $f\in L^\infty(X,\rho)$ define the operator $\pi(f)$ on
$L^2(X_\infty,\lambda_C)$ by
\begin{equation}\label{eqpi}
\pi(f)\xi(x,\omega)=f(x)\xi(x,\omega),\quad(\xi\in
L^2(X_\infty,\lambda_C), x\in X,\omega\in\Omega_x).
\end{equation}

\begin{proposition}\label{propcov}
The operator $U$ is unitary, $\pi$ is a representation of the
algebra $L^\infty(X,\rho)$ and the following covariance relation is
satisfied:
\begin{equation}\label{eqcov}
U\pi(f)U^{-1}=\pi(f\circ r),\quad(f\in L^\infty(X,\rho)).
\end{equation}
\par
For any function $f\in L^\infty(X,\rho)$ and $n\geq1$, the operator
$U^{-n}\pi(f)U^n$ is the operator of multiplication by
$f\circ\theta_n$. The union of the algebras $\{U^{-n}\pi(f)U^n\,|\,
f\in L^\infty(X,\rho)\}$ is SOT-dense in the algebra
$L^\infty(X_\infty,\lambda_C)$ (seen as multiplication operators on
$L^2(X_\infty,\lambda_C)$). An operator $S$ on
$L^2(X_\infty,\lambda_C)$ commutes with $U$ and $\pi(f)$ for all
$f\in L^\infty(X,\rho)$ if and only if there exists a function $F\in
L^\infty(X_\infty,\lambda_C)$ such that $F=F\circ\hat r$ and $S$ is
the operator of multiplication by $F$.
\end{proposition}

\begin{proof}
The fact that $U$ is an isometry follows form Proposition
\ref{proplambda}.
\par
The inverse of $U$ is
$$U^{-1}\xi(x,\omega)=\alpha_{i(\hat r(\omega))}^{-1}\frac{1}{\sqrt{\#r^{-1}(x)}}\xi\circ\hat r^{-1}(x,\omega).$$
Some simple computations prove the other relations.
\par
Note that the algebra $\{U^{-n}\pi(f)U^n\,\, f\in
L^\infty(X,\rho)\}$ is the algebra of operators of multiplication by
functions which depend only on the first $n$ coordinates. Since any
function in $L^\infty(X,\rho)$ can be pointwise and uniformly
boundedly approximated by such functions, it  follows that the union
of these algebras is dense in $L^\infty(X_\infty,\lambda_C)$.
\par
Since $L^\infty(X_\infty,\lambda_C)$ is maximal abelian, if $S$
commutes with $U$ and $\pi$ then $S$ commutes with the
multiplication operators, so it is a multiplication operator itself,
$S=M_F$. Since $S$ commutes with $U$, it follows that $F=F\circ\hat
r$.
\end{proof}
\par
Our formula for the measure $\lambda_C$ in (\ref{eqlambdac}) shows
that the Hilbert space $L^2(X_\infty,\lambda_C)$ fibers over
functions on $X$ as follows: for a dense space of functions
$\xi,\eta\in L^2(X_\infty,\lambda_C)$, the sum
$$\ip{\xi}{\eta}(x):=\sum_{\omega\in\mathbf{N}_C(x)}\xi(\omega)\cj{\eta}(\omega)$$
defines a $C(X)$-valued inner product as in \cite{Pac05},
\cite{HLPS99} and
$$\int_{X}\ip{\xi}{\eta}(x)\,d\rho(x)=\ip{\xi}{\eta}_{L^2(X_\infty,\lambda_C)}.$$
\subsection{A direct integral decomposition}\label{subsdire}
\par
We now resume our discussion of cycles $C$ for the endomorphism.
\par
 The cycle $C=\{x_0,\dots,x_{p-1}\}$ generates $p$ points in
$X_\infty$:
$$\omega_C:=(x_0,x_1,\dots,x_{p-1},x_0,\dots)\mbox{ and }\hat
r^k(\omega_C), k\in\{1,\dots,p-1\},$$ i.e., $\omega_C$ is the path
that goes through the cycle infinitely many times. We may write
$\omega_C:=CC\cdots=C^\infty$.

\begin{definition}\label{defrepel}
A fixed point $x_0$ for $r$ is called repelling if there is $0<c<1$
and $\delta>0$ such that for all $x\in X$ with $d(x,x_0)<\delta$,
one has $d(r(x),x_0)>c^{-1}d(x,x_0)$.\par A cycle
$C=\{x_0,\dots,x_{p-1}\}$ is called repelling if each point $x_i$ is
repelling for $r^p$.
\end{definition}

\begin{definition}\label{deffndsom}
A subset $A$ of $X_\infty$ is called a {\it cross section} if for
every path $\omega\in\cup_{x\in X}\mathbf{N}_C(x)\setminus\{\hat
r^k(\omega_C)\,|\,k\in\{0,\dots,p-1\}\}$, the intersection
$A\cap\{\hat r^k(\omega)\,k\in\bz\}$ contains exactly one point.
\end{definition}

\begin{proposition}
If $C$ is repelling, and $r$ is continuous at the points in $C$,
then there exists a cross section.
\end{proposition}

\begin{proof}
Using the continuity of $r$ and the fact that the cycle is
repelling, we can find a small $\delta>0$ and $0<c<1$ such that
$r^i(B(x_0,\delta))\cap B(x_0,\delta)=\emptyset$, for
$i\in\{1,\dots,p-1\}$, $r^{-p}(x_0)\cap B(x_0,\delta)=\{x_0\}$, and
such that $d(r^p(x),x_0)\geq c^{-1} d(x,x_0)$ for $x\in
B(x,\delta)$.
\par
Define
$$A:=\{(z_k)_{k\in\bz}\in X_\infty\,|\, z_0\in
r^p(B(x_0,\delta))\setminus B(x_0,\delta),\mbox{ and }z_{kp}\in
B(x_0,\delta)\mbox{ for }k\geq1\}.$$
\par
It is enough to prove that, for every path $(z_k)_{k\in\bz}$ in
$\mathbf{N}_C(x)$, except the special ones $\omega_C$ and the
others, there is a unique $k_0\in\bz$ such that
\begin{equation}\label{eqk0}
z_{k_0}\in r^p(B(x_0,\delta))\setminus B(x_0,\delta),\mbox{ and
}z_{kp}\in B(x_0,\delta)\mbox{ for }k\geq1.
\end{equation}
\par
Since $\omega$ is in $\mathbf{N}_C(x)$, the sequence $\{z_{kp}\}$
converges to one of the points $x_i$. Then, using the continuity of
$r$, $\{z_{kp+p-i}\}$ converges to $x_0$.
\par
Take the first $k_0\in\bz$ such that $z_{kp+k_0}\in B(x_0,\delta)$,
for all $k\geq1$. We still have to justify why there is a first one.
\par
If not, then $z_{-kp+k_0}\in B(x_0,\delta)$ for all $k\geq0$. Then,
using the fact that $x_0$ is repelling for $r^p$, there is a $c$
such that $0<c<1$, and for all $k\geq0$
$$\delta>d(z_{-kp+k_0},x_0)=d(r^{kp}(z_{k_0}),x_0)\geq c^{-k}
d(z_{k_0},x_0).$$ But then let $k\rightarrow\infty$, and obtain that
$z_{k_0}=x_0$. So, $z_{k_0-l}=r^l(x_0)=x_{l\mod p}$ for all
$l\geq0$. Also, since $z_{k_0+p}\in B(x_0,\delta)\cap r^{-p}(x_0)$,
we get $z_{k_0+p}=x_0$. By induction we obtain then that $\omega$ is
one of the special points in the orbit of $\omega_C$, which yields
the contradiction.
\par
Since $z_{k_0+p}\in B(x_0,\delta)$, $z_{k_0}=r^p(z_{k_0+p})$ is in
$r^p(B(x_0,\delta))$ but not in $B(x_0,\delta)$.
\par
Since $z_{k_0+kp}\in B(x_0,\delta)$ for $k\geq 1$, this proves that
$k_0+kp$ does not have the property (\ref{eqk0}).
\par
Since $z_{k_0}\not\in B(x_0,\delta)$, this proves that $k_0-kp$ does
not have the property (\ref{eqk0}) for $k\geq1$.
\par
Since for $k\geq 0$, $z_{k_0+kp+p}\in B(x_0,\delta)$, it follows
that for $i\in\{1,\dots,p-1\}$, one has $z_{k_0+kp+i}\in
r^{p-i}(B(x_0,\delta))$ so it is not in $B(x_0,\delta)$, and
therefore $k$ does not satisfy (\ref{eqk0}) when $k\not\equiv
k_0\mod p$. \par This proves that $A$ is a cross section.
\end{proof}
\par
Our present notion of cross section, and our next theorem are
motivated in part by an earlier theorem of Lim, Packer, and Taylor
\cite{LPT01} on direct integral decompositions of a class of wavelet
representations: This is the class of wavelets for which the Fourier
transform $\hat\psi$ of the wavelet mother-function $\psi$ is the
indicator function of a measurable subset in $\br^d$. Both our
present direct integral decomposition theorem, and that in
\cite{LPT01} are motivated by Mackey's theory of unitary
representations of non-abelian groups. In fact, our representation
of the covariant system $(U, \pi)$ may be viewed as a single
representation of a certain non-abelian crossed product
$\mathfrak{A}_{\hat r}:=C(X_\infty)\ltimes_{\hat r}\bz$ (see
\cite{BrJo91}), and our simultaneous direct integral decomposition
of $U$ and $\pi$ in Theorem \ref{thdecomp} below, is also a direct
integral decomposition of a single representation of the crossed
product group.

Assume now that $A$ is a cross section. For each $\omega\in A$,
define the operators $U_\omega$ and $\pi_\omega(f)$, $f\in
L^\infty(X,\rho)$ on $l^2(\bz)$ by
$$U_\omega\xi(k)=\alpha_{i(r^{k}(\omega))}\xi(k+1),\quad(\xi\in
l^2(\bz),k\in\bz),$$
$$\pi_\omega(f)\xi(k)=f(z_{-k})\xi(k),\quad(\xi\in
l^2(\bz),k\in\bz).$$
\par
The representation $\pi_\omega$ extends to a representation of
$C(X_\infty)$ defined by
$$\pi_\omega(f)\xi(k)=f(\hat r^k(\omega))\xi(k),\quad(f\in
C(X_\infty), \xi\in l^2(\bz),k\in\bz).$$ The covariance relation is
satisfied: $$U_\omega\pi_\omega(f)
U_\omega^{-1}=\pi_\omega(f\circ\hat r),\quad(f\in C(X_\infty)).$$
\begin{theorem}\label{thdecomp}
Let $A$ be a cross section, and assume that $\rho(C)=0$. The map
$\Phi\colon L^2(X_\infty,\lambda_C)\rightarrow L^2(A,\lambda_C)\otimes
l^2(\bz)$ defined by
$$(\Phi(f))(\omega,k)=\sqrt{\cc^{(k)}(z_0)}f(\hat
r^k(\omega)),\quad(f\in
L^2(X_\infty,\lambda_C),\omega=(z_k)_{k\in\bz}\in A, k\in\bz),$$ is
an isometric isomorphism which intertwines the operators $U$ and
$\int_A^\oplus U_\omega\,d\lambda_C(\omega)$, and also the
representations $\pi$ and
$\int_A^\oplus\pi_\omega\,d\lambda_C(\omega)$. The representations
$(U_\omega,\pi_\omega)$ on $l^2(\bz)$ are irreducible for all
$\omega\in A$.
\end{theorem}

\begin{proof}
The fact that $\Phi$ is isometric follows from Proposition
\ref{proplambda}. The inverse of $\Phi$ is defined as follows: for
each $\omega\in\cup_x\mathbf{N}_C(x)$, (except the special ones
which have measure $0$, so do not matter), there exists a unique
$k(\omega)\in\bz$ and $\eta(\omega)\in A$ such that $\omega=\hat
r^{k(\omega)}(\eta(\omega))$. Then
$$\Phi^{-1}(f)(\omega)=\frac{1}{\sqrt{\mathfrak{c}^{(k(\omega))}(\eta(\omega)_0)}}f(\eta(\omega),k(\omega)).$$
Everything follows by direct computation.
\par
We prove now that the representation $(U_\omega,\pi_\omega)$ is
irreducible for all $\omega=(z_k)_{k\in\bz}\in A$.
\par
Note first that $\pi_\omega(f)$ is a diagonal operator $F$ with
entries $F_{kk}=f(z_{-k})$, $k\in\bz$.
\par
We claim that for $k\neq l$ big enough, we have $z_k\neq z_l$. Since
$\omega$ is in $\mathbf{N}_C(z_0)$, it follows that $z_{kp}$
converges to one of the points of the cycle. Also, for $k$ big
enough, the points $z_{k}$ cannot be in $C$, because, this path
$\omega_C$ was removed from $A$. Suppose now that for any $m$ we can
find $k,l\geq m$, such that $k>l$ and $z_k=z_l$. Then this implies
that $z_k$ is periodic, therefore also
$z_{k-1}=r(z_k),z_{k-2},\dots,z_l$ are periodic, and since $m$ is
arbitrary, it follows that all the points $z_m$ are periodic. The
orbit of the two periodic points $z_0$ and $z_k$ intersect, because
$r^k(z_k)=z_0$, therefore the two orbits must be the same. Thus the
path $(z_k)_{k\in\bz}$ is an infinite repetition of the periodic
orbit of $z_0$: $(z_0,z_1,\dots,z_{p_0-1},z_0,z_1,\dots)$. However, this
cannot converge to the cycle $C$.
\par
Take now $k\neq l$ small enough. Then $z_{-k}\neq z_{-l}$ so we can
pick a function continuous function $f$ such that
$F_{kk}=f(z_{-k})\neq f(z_{-l})=F_{ll}$. If $T=(T_{ij})_{i,j\in\bz}$
is an operator on $l^2(\bz)$ that commutes with $U_\omega$ and
$\pi_{\omega}$, then $T_{kl}F_{ll}=F_{kk}T_{kl}$. So $T_{kl}=0$ for
$k\neq l$ small enough.
\par
Note that
$(U_\omega^{-1}\pi_\omega(f)U_\omega\xi)(k)=f(z_{-k+1})\xi(k)$, so
the conjugation with $U_\omega$ shifts the diagonal entries of
$\pi_\omega(f)$. Therefore, with the previous argument, we obtain
that $T_{kl}=0$ for all $k\neq l$. So $T$ is a diagonal operator.
Since, $T$ commutes with $U_\omega$, we obtain that
$T_{kk}=T_{k+1,k+1}$. So $T$ is a constant multiple of the identity.
This proves that the representation $(U_\omega,\pi_\omega)$ is
irreducible.
\end{proof}
\par
We show in Theorem \ref{thfaithful} below that the harmonic analysis
of covariant systems $(U,\pi)$ as in (\ref{eqcov}) is completely
equivalent to that of single representations $\hat\pi$ of a certain
$C^*$-algebraic crossed product $\mathfrak{A}_{\hat r}$. With this
identification $(U,\pi)\leftrightarrow\hat\pi$, we note in
particular that the operators in the commutant of the pair $(U,\pi)$
coincide with the commutant of the representation $\hat\pi$. Our
main conclusion in Theorem \ref{thfaithful} is that the
representation $\hat\pi$ of $\mathfrak{A}_{\hat r}$ is faithful,
i.e., that the kernel of $\hat\pi$ is trivial.
\begin{theorem}\label{thfaithful}
Assume that for every $x\in X$, there exists a path $(z_i)_{i\geq1}$
that starts at $x$ and with $\lim_{i\rightarrow\infty}z_{pi}\in C$,
(i.e., $\mathbf{N}_C(x)$ is non-empty). Then the operators $U$ and
$M_f$, $(f\in C(X_\infty))$ on $L^2(X_\infty,\lambda_C)$ form a
faithful representation of the crossed-product $\mathfrak{A}_{\hat
r}:=C(X_\infty)\ltimes_{\hat r}\bz$.
\end{theorem}
\begin{proof}
\par
The $C^*$-algebraic crossed product $\mathfrak{A}_{\hat r}$
\cite{Ped79} is the $C^*$-algebraic completion of formal symbols
$\{(f,k)\,|\,f\in C(X_\infty),k\in\bz\}$ with product
$$(f,k)\cdot(g,l)=(f\,g\circ\hat r^{k},k+l),\quad(f,g\in
C(X_\infty),k,l\in\bz).$$
\par
The representation is defined by
$$\hat\pi\colon(f,k)\mapsto \pi(f)U^k.$$
\par
We saw in Proposition \ref{propcov} and its proof that the
covariance relation is satisfied, so we have to check only that this
representation is faithful. If not, using a result from
\cite{Tak03}, we see that there is a non-zero element in
$\mathfrak{A}_{\hat r}$ of the form $(\sum_{k\in\bz}c_k(f,k))$ with
$\sum_{k\in\bz}|c_k|<\infty$, $f\in C(X_\infty)$ such that this
element is mapped to $0$ under $\hat\pi$.
\par
This means that $\pi(f)\sum_{k\in\bz}c_kU^k=0$. With Theorem
\ref{thdecomp} it follows that for almost all $\omega\in A$ and all
$\xi\in l^2(\bz)$, $l\in\bz$, one has
$$f(\hat r^{l}(\omega))\sum_{k\in\bz}c_k\xi(k+l)=0.$$
Take $\xi=\delta_i$ and get $f(\hat r^l(\omega))c_{i-l}=0$ which
implies that either $f(\hat r^l(\omega))=0$ for all $l$, or $c_l=0$
for all $l$. But, if $c_l=0$ for all $l$ then this contradict that
the element in the crossed-product in non-zero.
\par
Thus, for almost all $\omega=(z_i)_{i\in\bz}\in A$, we have that
$f(\hat r^l(\omega))=0$ for all $l$.
\par
This implies that $f$ is $0$ on almost all $\cup_x\mathbf{N}_C(x)$.
We know that non-empty open sets in $X$ have positive $\rho$-measure
(see \cite{DuJo05}). Hence, since the measure on each
$\mathbf{N}_C(x)$ is atomic, every non-empty open set in
$\cup_x\mathbf{N}_C(x)$ has positive measure. This implies that $f$
has to be $0$ on all $\cup_x\mathbf{N}_C(x)$ . We claim that this
set is dense in $X_\infty$.
\par
Take $\omega:=(z_1,z_2,\dots)\in X_\infty$ and $n\geq 1$ fixed. Since
$\mathbf{N}_C(z_n)$ is not empty, there exists a path
$(y_1,y_2,\dots)$ that starts at $z_n$ and is convergent to the cycle.
Then, $(z_1,z_2,\dots,z_n,y_1,y_2,\dots)$ is in $\mathbf{N}_C(x)$ and
coincides with the initial path on the first $n$ components. Thus,
the path $\omega$ can be approximated with paths in
$\mathbf{N}_C(x)$.
\par
Hence $\cup_x\mathbf{N}_C(x)$ is dense in $X_\infty$, and this
implies that $f=0$. The contradiction yields the result.

\end{proof}

\begin{remark}
{\bf[Iteration of rational functions]}
\par
Let $r\colon\mathbb{S}^2\rightarrow\mathbb{S}^2$ be a rational function
viewed as an endomorphism of the Riemann sphere $\mathbb{S}^2$, or
$\bc^\infty=\bc\cup\{\infty\}$; and suppose the degree of $r$ is
bigger than $2$. Let $X=X(r)$  be the {\it Julia set} of $r$, i.e.,
$X$ is the complement of the largest open subset $\mathcal{U}$ such
that $\{r^n|_{\mathcal{U}}\,|\,n\geq1\}$ is a normal family. It is
known that $X$ is non-empty, compact, and that $(X,r)$ carries a
unique strongly invariant probability measure; see \cite{Bea94} and
\cite{Bro65}.
\par
Our present general result for cycles are motivated by the following
specific theorems for rational mappings:
\par
Let $r$ be a rational mapping of degree at least $2$.
\begin{itemize}
\item Let $C$ be a $p$-cycle for $r$. Then $C$ is repelling if and
only if $|(r^p)'(z)|>1$ for all $z\in C$. Moreover,
$(r^p)'(z)=\prod_{w\in C}r'(w)$, $z\in C$, so $(r^p)'$ has the same
value for all points $z$ on the cycle $C$.
\item Every repelling cycle $C$ lies in the Julia set $X$.
\item The Julia set $X$ is the closure of the repelling periodic
points, see \cite[Theorem 3.1]{McMu94}.
\end{itemize}
\end{remark}

\subsection{The scaling function}\label{subsscal}
\par
We now turn to a theorem which is analogous to the existence theorem
for the scaling function in the classical theory of wavelets. As
outlined in \cite{Dau92}, the wavelet scaling function $\varphi$ in
$L^2(\br)$ depends on a filter function $m_0$ with $m_0$ defined on
$\bt=\br/\bz$. In fact, in the wavelet theory, it is the Fourier
transform $\hat\varphi$ which is an infinite product of scaled
versions of $m_0$. As is well known, this representation requires
that the function $m_0$ satisfies two conditions: one is called the
{\it quadrature condition}, and the second is called the {\it
low-pass} condition. Both of these conditions are motivated directly
from the probabilistic interpretation that $|m_0|^2$ enjoys in
signal processing.
\par
In our theorem below we identify the analogue of these two
conditions for the function $m_0\colon X\rightarrow\bc$ which is
associated to an endomorphism $r\colon X\rightarrow X$. The quadrature
condition takes the form (\ref{eqqmf}), and the low-pass condition
takes the form (\ref{eqlowpass}). The reason for the name quadrature
is that $r(z)=z^2$ in the wavelet case, and the reason for the name
low-pass, is that points on $\bt=\br/\bz$ correspond to frequencies,
and $x=0$ is the lowest frequency.
\par
In the general setting of the endomorphism $r$, the analogue of low
frequencies are points in cycles $C$ for $r$, and in this setting
low-pass means that $|m_0|^2$ attains its maximum on $C$. This is
exactly what condition (\ref{eqlowpass}) is saying.
\par
In the case of endomorphism, we will therefore expect to represent a
scaling function as an infinite product built out of $m_0$ and
iterated shifts applied to $m_0$. The fact that this can be done is
the main conclusion in the theorem.
\begin{definition}
A complex valued function $f$ on a metric space $X$ is called
$\beta$-Lipschitz at a point $x_0$ if there is a non-decreasing
function $\beta\colon[0,\infty)\rightarrow[0,\infty)$ such that for all
$A>0$ and $c<1$,
$$\sum_{k=1}^\infty\beta(Ac^k)<\infty,$$
and
$$|f(x)-f(x_0)|\leq\beta(d(x,x_0)),$$
for all $x$ in some neighborhood of $x$.
\end{definition}

\begin{definition}\label{def2_10}
Let $W\colon X\rightarrow[0,1]$ be a given function, and set
$$R_Wf(x):=\sum_{r(y)=x}W(y)f(y),\quad(x\in X).$$
We say that $R_W$ is a {\it transfer operator}, or a {\it Ruelle
operator}. A function $h$ on $X$ is said to be {\it harmonic} with
respect to $R_W$ if $R_Wh=h$. A cycle $C$ for $r$ is said to be a
$W$-cycle if $W(x)=1$ for all $x\in C$.
\end{definition}

The operator in Definition \ref{def2_10} plays a role in many areas of mathematics
and physics. Some of its recent uses are highlighted in Ref.\ \cite{Rue89}, where it is
key to Ruelle's thermodynamical formalism.

\begin{lemma}\label{lemrelforpx}
There is a unique family of measures $P_x$ supported on $\Omega_x$,
$x\in X$, satisfying the following relation: for all measurable sets
$E\subset X_\infty$ and all $x\in X$
$$\sum_{r(y)=x}W(y)P_y(\Omega_y\cap\hat r^{-1}(E))=P_x(\Omega_x\cap
E).$$
\end{lemma}
\begin{proof}
It is enough to define $P_x$ on cylinder sets: for a fixed
$(a_1,a_2,\dots,a_n,\dots)\in\Omega_x$ and $n\geq2$,
$$E:=\{(z_1,z_2,\dots)\in\Omega_x\,|\,z_1=a_1,\dots,z_n=a_n\}$$
Then $P_x(E\cap\Omega_x)=\prod_{k=1}^nW(a_k)$. \par The extension of
$P_x$ to the sigma-algebra generated by the cylinder sets now
follows from Kolmogorov's theorem. See \cite{DuJo05} for more
details. \par
 Also, for $y\in r^{-1}(x)$, one has that $\hat
r^{-1}(E)\cap\Omega_y$ is empty unless $y=a_1$, in which case
$P_y(\hat r^{-1}(E)\cap\Omega_y)=\prod_{k=2}^nW(a_k)$. The lemma
follows.
\end{proof}
\begin{lemma}\label{lem2dzdz}
The function $h_C(x):=P_x(\mathbf{N}_C(x))$ is harmonic with respect
to $R_W$.
\end{lemma}
\begin{proof}
We have the following disjoint union
$\cup_{r(y)=x}\mathbf{N}_C(y)=\hat r^{-1}(\mathbf{N}_C(x))$. The
lemma follows then from Lemma \ref{lemrelforpx}. Indeed,
$$(R_Wh_C)(x)=\sum_{r(y)=x}W(y)h_C(y)=\sum_{r(y)=x}W(y)P_y(\Omega_y\cap\mathbf{N}_C(y))$$
$$=\sum_{r(y)=x}W(y)P_y(\Omega_y\cap\hat
r^{-1}\mathbf{N}_C(x))=P_x(\Omega_x\cap\mathbf{N}_C(x))=h_C(x).$$
\end{proof}
\begin{definition}
We call $h_C(x):=P_x(\mathbf{N}_C(x))$ the harmonic function
associated to the cycle $C$. See also \cite{DuJo05} for more
details.
\end{definition}
\par
In our next theorem, we prove that each repelling cycle $C$
generates a covariant operator system on the corresponding Hilbert
space $L^2(X_\infty,\lambda_C)$. Moreover, under two conditions on a
given filter function $m_0$, we show that the corresponding scaling
equation has a natural solution $\hat\varphi_C$ in
$L^2(X_\infty,\lambda_C)$.
\begin{theorem}\label{propphic}
Assume now that the cycle $C$ is repelling and the functions $r$ and
$x\mapsto\#r^{-1}(x)$ are continuous at the points in $C$. Let
$m_0\in L^\infty(X,\rho)$ be a function which is $\beta$-Lipschitz
at the points in $C$ and which satisfies the conditions

\begin{equation}\label{eqqmf}
\frac{1}{\#r^{-1}(x)}\sum_{y\in r^{-1}(x)}|m_0(y)|^2=1,\quad(x\in X),
\end{equation}
and
\begin{equation}\label{eqlowpass}
m_0(x_i)=\alpha_i\sqrt{\#r^{-1}(r(x_i))},\quad(i\in\{0,\dots,p-1\}).
\end{equation}
Define the function $\hat\varphi$ by
\begin{equation}\label{eqphic}
\hat\varphi(x,(z_k)_{k\geq1}):=\prod_{k=1}^\infty\frac{\alpha_{i(\omega)+k}^{-1}m_0(z_k)}{\sqrt{\#r^{-1}(r(z_k))}},\quad(x\in X,(z_k)_{k\geq1}\in\mathbf{N}_C(x)).
\end{equation}
Then $\hat\varphi$ is in $L^2(X_\infty,\lambda_C)$ and it satisfies
the following relation:
\begin{equation}\label{eqscaling}
U\hat\varphi=\pi(m_0)\hat\varphi.
\end{equation}
If $W(x):=|m_0(x)|^2/\#r^{-1}(r(x))$, then $C$ is a $W$-cycle, and
if $h_C$ is the harmonic function associated to this $W$-cycle, then
\begin{equation}\label{eqcorr}
\ip{\pi(f)\hat\varphi}{\hat\varphi}=\int_Xfh_C\,d\rho.
\end{equation}
\par
Set $V_0:=\cj{\{\pi(f)\hat\varphi\,|\,f\in
L^\infty(X,\rho)\}}^{L^2}$, and $V_n:=U^{-n}V_0$, for $n\in\bz$.
Then $V_n\subset V_{n+1}$,
$$\cj{\bigcup_{n\in\bz}V_n}=L^2(X_\infty,\lambda_C),\quad\bigcap_{n\in\bz}V_n=\{0\}.$$
\end{theorem}

\begin{proof}
First we check that the infinite product (\ref{eqphic}) is
convergent. Take $x\in X$, $\omega=(z_1,z_2,\dots)\in\mathbf{N}_C(x)$.
Then the sequence $\{z_{kp}\}$ converges to one of the points of the
cycle $C$, namely $x_{i(\omega)}$. Applying $r$, which is continuous
at these points, we obtain that, for all $l$, the sequence
$\{z_{kp+l}\}$ is convergent to $x_{i(\omega)+l}$.
\par
Now we use the fact that the cycle is repelling. For $k$ large
enough, the path $\omega$ enters the neighborhood where the cycle is
repelling (see Definition \ref{defrepel}). Therefore, there are
constants $0<c_l<1$, $0\leq m_l<\infty$ such that for $k$ large
enough $d(z_{kp+l},x_{i(\omega)+l})\leq c_l^k m_l$, for all
$l\in\{0,\dots,p-1\}$. Take $c=\max\{\sqrt[p]{c_l}\}\in(0,1)$ and
$M:=c^{-p}\max\{m_l\}$ . Then for $k$ large enough
$$d(z_k,x_{i(\omega)+k})\leq c^k M.$$
\par
Since the function $\#r^{-1}(r(\cdot))$ is continuous at the points
of the cycle, we get that for $k$ large,
$\#r^{-1}(r(z_k))=\#r^{-1}(r(x_{i(\omega)+k}))=:A_k\geq1$.
\par
Let $\beta$ be the function given by the $\beta$-Lipschitz condition
for $m_0$ at all the points of the cycle (Take the minimum of these
functions over all the points of the cycle). Using the condition
$(\ref{eqlowpass})$, we have:
$$\left|\frac{\alpha_{i(\omega)+k}^{-1}m_0(z_k)}{\sqrt{\#r^{-1}(r(z_k))}}-1\right|=\frac{1}{\sqrt{A_k}}\left|m_0(z_k)-m_0(x_{i(\omega)+k})\right|$$$$\leq \beta(d(z_k,x_{i(\omega)+k}))\leq \beta(c^kM).$$
This implies that the sum over the terms on the left-hand side of
this inequality is convergent, which in turn implies that the
infinite product is absolutely convergent.
\par
Next we check (\ref{eqcorr}). It is clear that $C$ is a $W$-cycle. Also note that
$$|\hat\varphi(x,(z_k)_{k\geq1})|^2=\prod_{k=1}^\infty W(z_k)=P_x(\{(z_k)_{k\geq1}\}).$$
(See \cite{DuJo05}). Then
\begin{equation}\label{eqhcphic}
h_C(x)=P_x(\mathbf{N}_C(x))=\sum_{\omega\in\mathbf{N}_C(x)}|\hat\varphi(x,\omega)|^2,
\end{equation}
and equation (\ref{eqcorr}) follows. Since $h_C\leq 1$, this also
implies that $\hat\varphi$ is in $L^2(X_\infty,\lambda_C)$.
\par
We check the equation (\ref{eqscaling}). For $\omega=(z_1,z_2,\dots)\in\mathbf{N}_C(x)$,
$$U\hat\varphi(x,\omega)=\sqrt{\#r^{-1}(r(x))}\alpha_{i(\omega)}\prod_{k=1}^\infty\frac{\alpha_{i(\hat r(\omega))+k}^{-1}m_0(z_{k-1})}{\sqrt{\#r^{-1}(r(z_{k-1}))}}=$$
$$\sqrt{\#r^{-1}(r(x))}\alpha_{i(\omega)}\frac{\alpha_{i(\omega)}^{-1}m_0(x)}{\sqrt{\#r^{-1}(r(x))}}\prod_{k=2}^\infty\frac{\alpha_{i(\omega)+k-1}^{-1}m_0(z_{k-1})}{\sqrt{\#r^{-1}(r(z_{k-1}))}}=m_0(x)\hat\varphi(x,\omega).$$
\par
The scaling equation (\ref{eqscaling}) and the covariance equation
(\ref{eqcov}) imply that $V_{-1}\subset V_0$. This implies that the
sequence of subspaces $\{V_n\}$ is increasing.
\par
To check that their union is dense, we note that the closure of this
union is invariant for $U$ and for $\pi(f)$, $f\in
L^\infty(X,\rho)$. Therefore the projection $P$ onto this space is
in the commutant $\{U,\pi\}'$. But, then, with Proposition
\ref{propcov}, we obtain a function $F$ in
$L^\infty(X_\infty,\lambda_C)$ such that $F=F\circ\hat r$ and
$P=M_F$. In particular, $F\hat\varphi=\hat\varphi$ and $F$ is the
characteristic function of some set $\mathcal{F}$ which is invariant
for $\hat r$. \par However, the previous argument shows that, if
$\omega=(z_1,z_2,\dots)\in\mathbf{N}_C(x)$ has $z_i$ close enough to
the cycle $C$, for all $i\geq 1$, then $\hat\varphi(x,\omega)$ is
close to $1$. Now, take
$\omega=(z_1,z_2,\dots)\in\mathbf{N}_C(x)\setminus\mathcal{F}$. Then
$\hat r^{-n}(\omega)$ is outside $\mathcal{F}$. Because
$\omega\in\mathbf{N}_C(x)$, for $n$ large enough, all the points
$z_{n+1},z_{n+2},\dots$ are close to the cycle, so $\hat\varphi(\hat
r^{-n}(\omega))$ is close to $1$. But $\hat\varphi(\hat
r^{-n}(\omega))=\hat\varphi(\hat
r^{-n}(\omega))\chi_{\mathcal{F}}(\hat r^{-n}(\omega))=0$, a
contradiction. It follows that $\mathcal{F}$ has complement of
measure $0$ so $P=M_F$ is the identity, and therefore the union of
the multiresolution subspaces is dense.
\par
It remains to check that the intersection $\cap V_n$ is trivial. We
use the following lemma:
\begin{lemma}\label{lemjv0}
Define $\mathcal{J}\colon L^2(X,h_C\,d\rho)\rightarrow V_0$ by
$$\mathcal{J}(f)=\pi(f)\hat\varphi,\quad(f\in L^\infty(X,\rho)).$$
Define the operator $S_0$ on $L^2(X,h_C\,d\rho)$ by
$S_0f=m_0\,f\circ r$. Then $\mathcal{J}$ is an isometric isomorphism
such that $U\mathcal{J}=\mathcal{J}S_0$.
\end{lemma}
The proof of the lemma requires just some simple computations. The
fact that $S_0$ is an isometry is proved in Theorem \ref{thwold}.
\par
With this lemma, the assertion follows from Theorem \ref{thwold}.
\end{proof}
\par
Let $(X,\mathfrak{B},\rho)$ be a measure space with $\rho$ some
fixed probability measure defined on the sigma-algebra
$\mathfrak{B}$ on $X$. Let $\pi$ be a representation of
$L^\infty(X,\mathfrak{B})$ on a Hilbert space $\mathcal{H}$, and
suppose that the measure $f\mapsto\ip{\pi(f)\psi}{\psi}$ is
absolutely continuous with respect to $\rho$ for all
$\psi\in\mathcal{H}$, i.e., there exists $h_\psi\in
L^1(X,\mathfrak{B})$ such that
$\ip{\pi(f)\psi}{\psi}=\int_Xfh_\psi\,d\rho$, $f\in
L^\infty(X,\mathfrak{B})$. \par By the spectral multiplicity theorem
(\cite{BMM99}, \cite{Hal61}, \cite{NaFo70}), there is a measurable
function $d\colon X\rightarrow\{1,2,\dots,\infty\}$ such that if
$X_k:=\{x\in X\,|\,d(x)\geq k\}$, then the spectral representation
of $\pi$ takes the form of an isometric isomorphism
$\mathcal{J}\colon\mathcal{H}\rightarrow\sum_{k\in\bn}^{\oplus}L^2(X_k,\mathfrak{B},\rho)$,
such that
$\mathcal{J}_k\pi(f)\psi=f\mathcal{J}_k\psi=M_f\mathcal{J}_k\psi$
for all $f\in L^\infty(X,\mathfrak{B})$ and all
$\psi\in\mathcal{H}$.
\par
We say that $d$ is the {\it multiplicity function} of the
representation $\pi$.
\begin{corollary}
Let $V_0\subset L^2(X_\infty,\lambda_C)$, be the subspace from Lemma
\ref{lemjv0} and Theorem \ref{propphic}, and let $\pi_n$, $n\in\bn$,
be the restriction of the representation $\pi$ of
$L^\infty(X,\mathfrak{B})$ to $U^{-n}V_0$. Then
$$d_{U^{-n}V_0}(x)=\#(r^{-n}(x)\cap\{z\in X\,|\,h_C(z)\neq0\}),\quad(x\in X).$$
\end{corollary}
\begin{proof}
Since $\ip{\pi(f)\hat\varphi}{\hat\varphi}=\int_Xfh_C\,d\rho$, it
follows that $d_{V_0}(x)=\chi_{\{z\in
X\,|\,h_C(z)\neq0\}}=:\chi_{E_C}$.
\par
From \cite{DuJo04}, we know that
$$d_{U^{-n}V_0}(x)=\sum_{r^n(y)=x}d_{V_0}(y)=\sum_{r^{n}(y)=x}\chi_{E_C}(y)=\#(r^{-n}(x)\cap
E_C).$$
\end{proof}
 \begin{example}
Let $A$ be a square matrix with $0-1$ entries. Suppose every column
of $A$ contains at least one entry $1$. Then we show that the two
systems $(X_\infty,\hat r)$ and $(X,r)$ may be realized as
two-sided, respectively one-sided subshifts.
\par
Let $I$ be the index set for the rows and columns of $A$. Let
$$X_\infty(A):=\{(x_i)_{i\in\bz}\in I^{\bz}\,|\,
A(x_i,x_{i+1})=1\}.$$ Let $$\hat
r((x_i)_{i\in\bz})=(x_{i+1})_{i\in\bz}.$$ Define
$\theta_0((x_i)_{i\in\bz})=(x_i)_{i\geq0}$, and set
$X(A)=\theta_0(X_\infty(A))$.
\par
Then there is an endomorphism $r=r_A\colon X(A)\rightarrow X(A)$ such that
$r\circ\theta_0=\theta_0\circ\hat r$.
\par
Specifically, $$x=(x_i)_{i\in\bz}=\dots x_{-2}x_{-1}x_0x_1x_2\dots$$ with
$x_i\in I$;
$$\theta_0((x_i)_{i\in\bz})=(x_i)_{i\geq0}=x_0x_1x_2\dots;$$ and
$r(x_0x_1x_2\dots)=(x_1x_2x_3\dots)$ for $x\in X(A)$.
\par
For $x,y\in X(A)$, let $x\wedge y$ be the longest initial block in
$I\times I\times\cdots$ common to both $x$ and $y$, and let $|x\wedge
y|$ be the length of this block. Let $0<c<1$, and set
$d_c(x,y)=c^{|x\wedge y|}$. Then $d_c$ is a metric, and $(X(A),d_c)$
is a compact metric space whose open sets are generated by the
cylinder sets in $X(A)$. Moreover, $d_c(r(x),r(y))\leq
c^{-1}d_c(x,y)$ holds for all $x,y\in X(A)$. If $x\in X(A)$, then
$r^{-1}(x)=\{(ix)\,|\, A(i,x_0)=1\}$, and for the transfer operator
$\mathcal{L}_A\colon C(X(A))\rightarrow C(X(A))$,
$$(\mathcal{L}_Af)(x)=\frac{1}{\#r^{-1}(x)}\sum_{r(y)=x}f(y),$$
we have
$$(\mathcal{L}_Af)(x)=\frac{1}{\#\{i\,|\,A(i,x_0)=1\}}\sum_{A(i,x_0)=1}f(ix).$$
By \cite[Chapter 2]{MaUr03}, there is a unique probability measure
$\rho=\rho_A$ on $X(A)$ such that $\rho(\mathcal{L}_Af)=\rho(f)$ for
all $f\in C(X(A))$; i.e., $\rho=\rho_A$ is the unique strongly
invariant probability measure on $X(A)$.
\par
It follows that all the results in this setting apply; in
particular, if $C\subset X(A)$ is a cycle, then
$L^2(X_\infty(A),\lambda_C)$ is defined by
$$\int_{X_\infty(A)}|f|^2\,d\lambda_C=\int_{X(A)}\sum_{\omega\in\mathbf{N}_C(x)}|f(\omega)|^2\,d\rho(x)<\infty.$$
\par
Note also that $\mathbf{N}_C(x)$ consists of doubly infinite words
in $X_\infty(A)$ that start with an infinite repetition of the cycle
$C$. Specifically, for $x=(x_0x_1x_2\dots)\in X(A)$,
$\mathbf{N}_C(x)=\left\{(\omega_i)_{i\in\bz}\in
X_\infty(A)\,|\,\exists k\in\bn\right.$ such that
$(\omega_I)_{i\leq-k}\mbox{ is }C^\infty$, $(\omega_i)_{-k<i\leq-1}$
is some finite word, and $\left.\omega_i=x_i\mbox{ for
}i\geq0\right\}$.
\par
We now turn to a concrete example. Let the index set $I$ be
$\{1,2\}$ and let
$A=\left(\begin{array}{cc}1&1\\1&0\end{array}\right)$. This is
called the golden mean shift (see \cite[page 37]{LiMa95}).

\end{example}
\begin{proposition}\label{pr2_18} Let $m_0$ be the function on $X(A)$ determined by
\begin{equation}\label{eqexgm}
m_0(11\dots)=\sqrt{2},\quad m_0(21\dots)=0,\quad m_0(12\dots)=1,
\end{equation}
and $C$ be the cycle $\{111\dots\}$. Then $m_0$ satisfies
(\ref{eqqmf}) and (\ref{eqlowpass}) and defines a scaling vector
$\hat\varphi\in L^2(X_\infty,\lambda_C)$ with $h_C=1$.
\end{proposition}
\begin{proof}It is easy to
verify the two conditions (\ref{eqqmf}) and (\ref{eqlowpass}). The
scaling function $\hat\varphi$ is defined by the infinite product
(\ref{eqphic}). (We set $\alpha_i=1$.) If $\omega$ is in
$\mathbf{N}_C(x_0x_1\dots)$, then it has the form
$\dots111x_{-n}x_{-n+1}\dots x_{-1}x_0x_1\dots$. Note that if one of the
letters $x_{-k}$ is $2$ ($k\geq1$) then the next one has to be $1$.
Therefore, shifting the word to the right will bring the $21$ to the
central position, and $m_0$ is $0$ on words that start with $21$.
Therefore the infinite product is non-zero only when $x_{-k}=1$ for
all $k\geq1$. Then, an analysis of the possibilities for $x_0$ shows
that $\hat\varphi$ is $1$ in this case.
\par
Therefore $\hat\varphi(\omega)=1$, if $x_{-k}=1$ for all $k\geq 1$,
and $\hat\varphi(\omega)=0$ otherwise.
\par
Then by (\ref{eqhcphic})
$$h_C(x_0x_1\dots)=\sum_{\omega\in\mathbf{N}_C(x_0x_1\dots)}|\hat\varphi(\omega)|^2=1.$$
\par
An interesting consequence of (\ref{eqqmf}) and (\ref{eqlowpass})
for this example is that an admissible $m_0$ cannot be of the form
$m_0=\sqrt{2}\chi_E$ for a subset $E$ of $X(A)$ (because
$|m_0(21\ldots)|^2=1$). This contrasts a known wavelet, the Shannon
wavelet, see \cite{HeWe96} and \cite{BJMP05}.
\end{proof}

\section{Ergodic properties and the Wold decomposition}
\par
In our analysis of the intersection of the multiresolution spaces $V_n$, we are forced to study some convergence properties for the measure $\rho$ and the
filter $m_0$. The main tool in this study will be Doob's convergence theorems for reversed martingales, see e.g. \cite{Nev75}.
\par
Even though we are mainly interested in the strongly invariant
measure $\rho$, our analysis works in the following more general
case.
\begin{definition}
Let $V\geq0$ be a measurable function on $X$ such that
$$\sum_{r(y)=x}V(y)=1,\quad(x\in X).$$
A probability measure $\nu$ on $X$ such that
\begin{equation}\label{eqnu}
\int_Xf\,d\nu=\int_X\sum_{r(y)=x}V(y)f(y)\,d\nu(x),\quad(f\in L^1(X,\nu)).
\end{equation}
is called a Perron-Frobenius measure for the corresponding Ruelle
operator $$(R_Vf)(x):=\sum_{r(y)=x}V(y)f(y),\quad(x\in X).$$ For
example, when $V(x)=1/\#r^{-1}(r(x))$, then (\ref{eqnu}) is
equivalent to the strong invariance of $\nu$.
\end{definition}
\par
Note that a Perron-Frobenius measure $\nu$ is also invariant for
$r$, because
$$\int_Xf\circ r\,d\nu=\int_X\sum_{r(y)=x}V(y)f(r(y))\,d\nu(x)=\int_Xf(x)\sum_{r(y)=x}V(y)\,d\nu(x)=\int_Xf\,d\nu.$$
\par
Let $\mathfrak{B}$ be the sigma-algebra of measurable subsets of $X$.
\begin{definition}
Let $\mathfrak{B}$ be a sigma-algebra on $X$ and let $\nu$ be a
probability measure on $(X,\mathfrak{B})$. Let
$\mathfrak{C}\subset\mathfrak{B}$ be a sub-sigma-algebra. Then the
$\mathfrak{C}$-conditional expectation $E_{\mathfrak{C}}$ is defined
by
$$\int_XE_{\mathfrak{C}}f\,g\,d\nu=\int_Xfg\,d\nu,$$
for $f\in L^1(\mathfrak{B},\nu)$, $g\in L^\infty(\mathfrak{C})$; and
$E_{\mathfrak{C}}L^1(\mathfrak{B},\nu)=L^1(\mathfrak{C},\nu)$.
\end{definition}
\begin{proposition}\label{propcondexp}
The operator $E_n^V$ defined on $L^1(X,\nu)$ by
$$E_n^V(f)(x)=\sum_{r^n(y)=r^n(x)}V^{(n)}(y)f(y),\quad(x\in X),$$
defines the conditional expectation of $\mathfrak{B}$ with respect
to $r^{-n}(\mathfrak{B})$.
\end{proposition}

\begin{proof}
First note that if a function $g$ on $X$ is
$r^{-n}(\mathfrak{B})$-measurable, then $g(x)=g(y)$ whenever
$r^n(x)=r^n(y)$. Take now $g\in L^2(r^{-n}(\mathfrak{B}))$ and $f\in
L^1(\mathfrak{B})$. Using the invariance of $\nu$ and (\ref{eqnu}),
we have
$$\int_XE_n^V(f)g\,d\nu=\int_X\sum_{r^n(x)=r^n(y)}V^{(n)}(y)f(y)g(x)\,d\nu(x)$$$$
=\int_X\sum_{r^n(y)=x}V^{(n)}(y)f(y)g(y)\,d\nu(x)=\int_Xfg\,d\nu.$$
This shows that $E_n^V$ is the conditional expectation.
\end{proof}
\par
We note the relation between the Ruelle operator $R_V$ and the
conditional expectation $E_n^V$:
\begin{equation}\label{eqrvev}
E_n^V(f)=(R_V^n)\circ r^n,\quad(n\geq1, f\in L^1(X,\nu)).
\end{equation}
\par
The sigma-algebras $r^{-n}(\mathfrak{B})$ form a decreasing sequence, and we denote their intersection by $\mathfrak{B}_\infty$. Denote by $E_\infty^V$ the conditional expectation of $\mathfrak{B}$ with respect to $\mathfrak{B}_\infty$. Doob's theorems for reverse martingales can be applied now directly and we obtain the following theorem:

\begin{theorem}\label{thdoobreverse}
If $f\in L^p(X,\nu)$, $(1\leq p<\infty)$, then $E_n^V(f)$ converges
pointwise $\nu$-a.e. and in $L^p(X,\rho)$ to $E_\infty^V(f)$.
\end{theorem}

\begin{definition}\label{defaver}
We say that $r$ is {\it averaging} (with respect to the measure
$\nu$), if $L^1(\mathfrak{B}_\infty)$ contains only functions which
are constant $\nu$-a.e., (or, equivalently, the sigma-algebra
$\mathfrak{B}_\infty$ contains only sets of $\nu$-measure $0$ or
$1$).
\end{definition}

\begin{proposition}\label{propavererg}
If $r$ is averaging with respect to $\nu$ then it is also ergodic with respect to $\nu$.
\end{proposition}
\begin{proof}
If $A$ is a completely invariant set for $r$ then, for any two points $x,y$ such that $r^n(x)=r^n(y)$ for some $n\geq0$ $\chi_A(x)=\chi_A(r^n(x))=\chi_A(r^n(y))=\chi_A(y)$, so $\chi_A\in L^1(\mathfrak{B}_\infty)$, therefore $\nu(A)$ is $0$ or $1$.
\end{proof}

\begin{corollary}\label{cordoobaver}
If $r$ is averaging with respect to $\nu$, then for all $f\in L^p(X,\nu)$, $1\leq p<\infty$, the sequence $E_n^V(f)$ converges pointwise $\nu$-a.e. and in $L^p(X,\nu)$ to $\int_Xf\,d\nu$.
\end{corollary}

\par
Next, we will derive an ergodic property for a function $m_0$ satisfying (\ref{eqqmf}).
\begin{theorem}\label{thergm0}
Assume that the strongly invariant measure $\rho$ is ergodic with respect to $r$. Let
$m_0\in L^\infty(X,\rho)$ be a function that satisfies (\ref{eqqmf})
and such that $|m_0|\neq 1$ on a set of positive measure. Then

$$A:=\int_X\ln|m_0(x)|\,d\rho(x)\in[-\infty,0).$$
Then
$$\lim_{n\rightarrow\infty}|m_0(x)\cdots m_0(r^{n-1}(x))|^{\frac1n}=e^A,\mbox{ for }\rho\mbox{-a.e. }x\in X.$$
\end{theorem}
\begin{proof}
We have, using the strong invariance of $\rho$:
$$\int_X\ln|m_0(x)|\,d\rho(x)=\frac12\int_X\ln|m_0(x)|^2\,d\rho(x)=\int_X\frac{1}{\#r^{-1}(x)}\sum_{r(y)=x}\ln|m_0(y)|^2\,d\rho(x)$$
$$=\int_X\ln\left(\prod_{r(y)=x}|m_0(y)|^2\right)^{\frac1{\#r^{-1}(x)}}\,d\rho(x)$$$$\leq\int_X\ln\left(\frac1{\#r^{-1}(x)}\sum_{r(y)=x}|m_0(y)|^2\right)\,d\rho(x)
=\int_X\ln(1)=0.$$
 If we have equality in this chain, then we get that for $\rho$-a.e., $x\in X$, $|m_0(y)|=|m_0(y')|$ for all $y,y'\in r^{-1}(x)$, which implies that
$$1=\frac{1}{\#r^{-1}(r(x))}\sum_{r(y)=r(x)}|m_0(y)|^2=|m_0(x)|,\mbox{ for a.e. }x.$$
This contradicts the hypothesis. Thus $A\in[-\infty,0)$.
\par
Assume now, that $A>-\infty$. Then, using Birkhoff's ergodic
theorem, we obtain that
$$\lim_{n\rightarrow\infty}\ln\left(\frac{|m_0(x)\cdots m_0(r^{n-1}(x))|^{1/n}}{e^{A}}\right)=\lim_{n\rightarrow\infty}\frac{1}{n}\sum_{k=0}^{n-1}\ln|m_0(r^k(x))|-A$$
$$=\int_X\ln|m_0(x)|\,d\rho(x)-A=0.$$
 This yields the conclusion in the case $A>-\infty$.
\par
When $A=-\infty$, take $0>B>-\infty$ arbitrary and choose a bounded
measurable function $f$, with $|f|\geq|m_0|$ and such that
$-\infty<\int_X\ln|f(x)|\,d\rho(x)=B$. Then apply the previous
argument to conclude that $|f(x)f(r(x))\cdots f(r^{n-1}(x))|^{1/n}$
converges a.e. to $e^B$. Then
$$\limsup_n|m_0(x)\cdots m_0(r^{n-1}(x))|^{1/n}\leq e^B$$
and, as $B$ is arbitrary this implies that the limit is
$e^{-\infty}=0$.
\end{proof}

\par
With these results, we are now able to derive the result about the
Wold decomposition \cite{NaFo70} of the isometry $S_0$ associated to
$m_0$.

\begin{theorem}\label{thwold}
Let $\rho$ be a strongly invariant measure for $r$. Let $m_0\in
L^\infty(X,\rho)$ be a function that satisfies (\ref{eqqmf}). Let
$h\in L^\infty(X,\rho)$ be a function such that $h\geq0$ and
\begin{equation}\label{eqqmfh}
\frac{1}{\#r^{-1}(x)}\sum_{r(y)=x}|m_0(y)|^2h(y)=h(x),\quad(x\in X).
\end{equation}
Then the operator $S_0$ on $L^2(X,h\,d\rho)$ defined by
$$S_0f=m_0f\circ r$$
is an isometry.
\par
Assume in addition that $r$ is averaging with respect to $\rho$, and
that $|m_0|\neq 1$ on a set of positive measure $\rho$.  Then
$$\cap_{k\geq 1}S_0^k(L^2(X,h\,d\rho))=\{0\}.$$
\end{theorem}
\begin{proof}
The fact that $S_0$ is an isometry follows from the fact that $\rho$
is strongly invariant and from the relation (\ref{eqqmfh}):
$$\int_X|m_0(x)|^2|f(r(x))|^2h(x)\,d\rho(x)=\int_X\frac{1}{\#r^{-1}(x)}\sum_{r(y)=x}|m_0(y)|^2|f(r(y))|^2h(y)\,d\rho(x)$$$$=\int_X|f(y)|^2\,d\rho(x).$$
Denote, by
$$c(x):=\frac{1}{\#r^{-1}(r(x))},\quad(x\in X).$$
Note that
$$R_{m_0}^kf(x)=\sum_{r^k(y)=x}c^{(k)}(y)|m_0^{(k)}(y)|^2f(y),$$
where $R_{m_0}$ is the Ruelle operator associated to
$W(x):=|m_0(x)|^2/\#r^{-1}(r(x))$.

In particular
$$\sum_{r^k(y)=x}c^{(k)}(y)|m_0^{(k)}(y)|^2=1.$$
\par
Take now $\xi\in\cap_kS_0^k(L^2(X,h\,d\rho))$. Then for all $k\geq
1$, there exists $f_k\in L^2(X,h\,d\rho)$ such that
$\xi=m_0^{(k)}f_k\circ r^k$. This implies that for all $x\in X$:
\begin{align*}
|\xi(x)|^2&=|m_0^{(k)}(x)|^2|f_k(r^k(x))|^2\sum_{r^{k}(y)=r^k(x)}c^{(k)}(y)|m_0^{(k)}(y)|^2\\
&=|m_0^{(k)}(x)|^2\sum_{r^k(y)=r^k(x)}c^{(k)}(y)|m_0^{(k)}(y)f_k(r^k(y))|^2\\
&=|m_0^{(k)}(x)|^2\sum_{r^k(y)=r^k(x)}c^{(k)}(y)|\xi(y)|^2\\
&=|m_0^{(k)}(x)|^2E_k^c(|\xi|^2).
\end{align*}
With Theorem \ref{thergm0} and Corollary \ref{cordoobaver}, if we
let $k\rightarrow\infty$, we can conclude that $\xi=0$, $\rho$-a.e.
This proves the theorem.
\end{proof}

\begin{remark}
It is conceivable that the last conclusion in Theorem \ref{thwold} above may hold slightly more generally; possibly when only ergodicity
is assumed for $(X, r, \rho)$. But for the applications we have in mind, our
present assumption of strong invariance is appropriate, i.e., the
averaging assumption we place on the system $(X, r, \rho)$.
\end{remark}

\subsection{Some conditions for $r$ to be averaging}
We will give some necessary conditions for $r$ to averaging. For
this we will relate the expectation $E_V^n$ to the Ruelle operator
$R_V$.
\par
Just as before, assume $V\geq0$ is a measurable function such that
$$\sum_{r(y)=x}V(y)=1,\quad(x\in X),$$
and let $\nu$ be a measure such that
$$\int_XR_Vf\,d\nu=\int_Xf\,d\nu.$$

\begin{proposition}\label{propcondaver}
Suppose there exists a family of functions $\mathcal{F}$ which is
dense in $L^1(X,\nu)$ such that for all $f\in\mathcal{F}$,
$$\lim_{n\rightarrow\infty}\|R_V^n(f)-\int_Xf\,d\nu\|_1=0.$$
Then, for all $f\in L^1(X,\nu)$.
$$\lim_{n\rightarrow\infty}R_V^n(f)=\int_Xf\,d\nu=E_\infty^V(f).$$
In particular $r$ is averaging with respect to $\nu$.
\end{proposition}
\begin{proof}
Take $f\in L^1(X,\nu)$, and $\epsilon>0$. There exists
$g\in\mathcal{F}$, such that $\|f-g\|_1<\epsilon$. Then, using the
fact that $\nu$ is invariant for $r$, and also for $R_V$, we have,
with the aid of (\ref{eqrvev}):
$$\|E_n^V(f)-\int_Xf\,d\nu\|_1=\|R_V^nf-\int_Xf\,d\nu\|_1$$$$\leq\|R_V^n(f-g)\|_1+\|R_V^ng-\int_Xg\,d\nu\|_1+\|\int_X(g-f)\,d\nu\|_1$$
$$\leq 2\|f-g\|_1+\|R_V^ng-\int_Xg\|_1<3\epsilon,$$
for $n$ large enough. This proves the first assertion. Since
$E_\infty^V(f)$ is constant for all $f\in L^1(X,\nu)$, it follows
that $L^1(\mathfrak{B}_\infty)$ contains only constant functions so
$r$ is averaging.
\end{proof}

\begin{remark}\label{remaver}
The conditions of Proposition \ref{propcondaver} are satisfied in
many cases of interest. This is a consequence of Ruelle's theorem
(see \cite{Ba00}, \cite{FaJi01}). For example, if $r$ is locally
expanding (i.e., there exists $b>0$ and $\lambda>1$ such that
$d(r(x),r(y))\geq\lambda d(x,y)$ when $d(x,y)<b$), and mixing (i.e.,
for every open set $U$ in $X$, there exists $n$ such that
$r^n(U)=X$), and if $V>0$ and is Lipschitz, then $\mathcal{F}$ can
be taken to be the set of continuous functions, and $R_V^nf$
converges uniformly to $\int_Xf\,d\nu$, where $\nu$ is the unique
probability measure invariant for $R_V$.
\par
In particular, this is satisfied, for subshifts of finite type.
\par
Also, consider the case when $r$ is a rational map on $\bc$ and $X$
is its Julia set. Take $V=1/N$ where $N$ is the degree of the map
$r$. Then $\nu=\rho$ is the unique strongly invariant measure and we
may take again $\mathcal{F}$ to be the set of continuous functions
(see \cite{Lj83}). \par Given our assumptions above, the existence
and the uniqueness of the measure $\nu$ follows from the conclusion
in Ruelle's theorem, applied to $R_V$.

\end{remark}
\begin{acknowledgements} We thank Professors David Kribs, Kathy Merrill, Judy Packer, and Paul
Muhly for helpful discussions. Professor Paul Muhly enlightened us about a
number of places in the literature where variants of the general problem of
extension from non-invertible dynamics to a bigger ambient system occur. In
different contexts, this comes up in for example the papers \cite{BBLS04}, \cite{KaKr05} and
\cite{MuSo03}; as well as in the papers cited there.
\end{acknowledgements}

\end{document}